\documentclass[a4paper,12pt,dvipdfmx]{article}
\usepackage{graphicx}

\usepackage{amsmath,amssymb,cite,comment}

\usepackage{amsfonts}

\def\beq{\begin{equation}}
\def\eqn#1{\beq\label{#1}}
\def\ee{\end{equation}}
\def\eeq{\end{equation}}

\def\bb {\begin {eqnarray}}
\def\eqnn#1{\bb\label{#1}}
\def\eea {\end {eqnarray}}

\newcommand{\eqna}[1]{\begin{subequations} \label{#1}
\begin{eqnarray}}
\def\eena{\end{eqnarray}
\end{subequations}}

\def\th{{\tilde h}}

\def\nn{\nonumber}

\def\nd{\end{document}}

\input epsf.tex
\newcount\figno
\def\fig#1#2#3{
\par\begingroup\parindent=0pt\leftskip=1cm\rightskip=1cm\parindent=0pt
\baselineskip=11pt \global\advance\figno by 1 
\epsfxsize=#3 \centerline{\epsfbox{#2}} \vskip 12pt
#1\par
\endgroup\par}
\def\figlabel#1{\xdef#1{\the\figno}}
\def\encadremath#1{\vbox{\hrule\hbox{\vrule\kern8pt\vbox{\kern8pt
\hbox{$\displaystyle #1$}\kern8pt} \kern8pt\vrule}\hrule}}

  \def\tV{{\tilde V}}

\def\rank{{\rm rank}}

\def\riga{-\kern-4pt - \kern-4pt -}
\font\fat=cmsy10 scaled\magstep5

\def\Bbullet{\raise-3pt\hbox{\fat\char"0F}}

\font\tfont=cmbx12 scaled\magstep1 

\def\Box{
\vbox{ \halign to5pt{\strut##& \hfil ## \hfil \cr &$\kern -0.5pt
\sqcap$ \cr \noalign{\kern -5pt \hrule} }}~}

\def\down{\raise1.5pt\hbox{$\phantom{a}_2$}\downarrow}

\def\downa{\raise1.5pt\hbox{$\phantom{a}_{2\atop m_2}$}\downarrow}

\def\llr{\longrightarrow}

\def\({\left(}
\def\){\right)}

\def\ha{{\textstyle{1\over2}}}

\def\bbc{\mathbb{C}}
\def\bac{\mathbb{C}}

\def\bbr{\mathbb{R}}
\def\bbn{\mathbb{N}}

\def\a{\alpha}
\def\b{\beta}
\def\d{\delta}

\def\vr{\vert}

\def\l{\lambda}

\def\D{{\Delta}}

\def\ca{{\cal A}}  \def\cc{{\cal C}}
\def\cd{{\cal D}} \def\ce{{\cal E}} \def\cf{{\cal F}}
\def\cg{{\cal G}} \def\ch{{\cal H}} 
 \def\ck{{\cal K}} 
\def\cm{{\cal M}} \def\cn{{\cal N}} 
\def\cp{{\cal P}}  
 \def\ct{{\cal T}}

\def\ido{intertwining differential operator}
\def\idos{intertwining differential operators}

\def\L{\Lambda}
\def\r{\rho}

\begin{document}

\begin{center}

{\tfont Langlands Duality and Invariant Differential Operators}

\vskip 1.5cm

{\bf V.K. Dobrev}

 \vskip 5mm

  Institute for Nuclear Research and Nuclear Energy,\\ Bulgarian
Academy of Sciences,\\ 72 Tsarigradsko Chaussee,  1784 Sofia,
Bulgaria

\end{center}

\vskip 1.5cm

 \centerline{{\bf Abstract}}

Langlands duality is one of the most influential topics in mathematical research.  It has many different appearances and influential subtopics. Yet there is a topic that until now seems unrelated to the Langlands program.  That is the topic of invariant differential operators. That is strange since both items are deeply rooted in Harish-Chandra's representation theory of semisimple Lie groups.  In this paper we start building the bridge between the two programs.

\vskip 1.5cm

\section{Introduction}

In the last 50 years Langlands duality is one of the most influential topics in mathematical research \cite{Lan1,Lan2}.  It has many different appearances and influential subtopics, cf. an incomplete list 	\cite{AdVo,AFO,ABHL,ArFr,BaTe,Beil,BeNa,BSV,BrKa,CaRa,ChNa,ChVe,Chua,DHKM,DiTe,DoPa,Dri,Esp,EFK,Far,FaSch,FeFr,Fre,FrGa,FGKV,FGV,FrHe,GaTe,GaWi,Gin,GuWi,HLM,Han,Haus,Hit,Hove,Ike,Ima,JLN,KaWi,KiNo,KSZ,Laf,LZZ,MaOp,Mat,Nak,Suz,Tai,Tan,Tes,Wit}. Note that some papers are written by authors who have created influential topics themselves. The last fact
stresses the omnipresence of the Langlands program.

Yet there is a topic that until now seems unrelated to the Langlands program.  That is the topic of ~{\it invariant differential operators}. That is strange since both items are deeply rooted in Harish-Chandra's representation theory of semisimple Lie groups.  In this paper we start building a bridge between the two programs.

Our attempt is based on our approach to the construction of invariant differential operators - for an exposition we refer to \cite{VKD1} which is based also on many papers, see loc.cit.
Our approach is deeply related to the Langlands general classification of  representations of real semisimple groups $G$ \cite{Lan2} taking into account the refinement by Knapp-Zuckermann \cite{KnZu}.

One main ingredient in \cite{Lan2} are the parabolic subgroups $P=MAN$, such that $M$ is semisimple subgroup of our group $G$ under study, $A$ is abelian subgroup, $N$ is nilpotent subgroup preserved by the action $A$. Altogether, there is a local (Bruhat) decomposition of $G$ using   a subgroup $G'= P\tilde{N}$, where $\tilde{N}$ is a  nilpotent subgroup of $G$ isomorphic to $N$ also preserved by the action $A$, so that $G'$ is dense in $G$.
According to the construction of Langlands-Knapp-Zuckermann every admissible irreducible representation  of $G$ may be obtained as a subrepresentation of representations of $G$ induced by a representations of some $P$ (some class is enough - see details below).

Our construction of \idos\ is based on the fact that the structure of parabolic subgroups is related to various subgroups of the Weyl groups   $W(\cg^\bbc,\ch^\bbc)$, where $\cg$ is the Lie algebra of $G$,
$\ch$ is the Cartan subalgebra of some $MA$. This is also related to
various intertwining operators  in the Langlands dual group, cf. \cite{EFK,BrKa}.

 Another aspect of the above is the Chevalley automorphism  in the case of real groups which is an exhibition   of the local Langlands correspondence over R \cite{AdVo}. This is also related to the   notion of Hermitian dual, see the case of $SL(n,\bbr)$ below. Another application to representation theory using Heisenberg modules is done in \cite{Beil}.

 Gauge theory aspects  of the geometric Langlands program are studied in \cite{GaWi,GuWi,KaWi,Tai,Wit}.
 An exotic  'Chtoucas' application of the Langlands program is given in \cite{Laf}.
 Similarly, Langlands duality extends to Poisson-Lie duality via cluster theory. \cite{ABHL}. Further,   Langlands duality extends of representations of W-algebras in the quantum framework \cite{ArFr}.
 A proof of the global Langlands conjecture for GL(2) over a function field is given in \cite{Dri}.

  Further a two-parameter generalization of the geometric Langlands correspondence is   proved for all simply-laced Lie algebras \cite{AFO}. This is related to two-parameter quantum groups, e.g., \cite{VKD2}, and to
   6d conformal   supersymmetry \cite{VD6}. The conformal case is studied also in \cite{FrGa}.
   The supersymmetry case is studied also in \cite{BaTe}.

Various aspects of the Langlands program are given in \cite{Fre}.
The local Langlands correspondence is studied in \cite{Far,FaSch}.
Applications to integrability are shown in \cite{DoPa,FeFr,FrHe,Tan,Tes}.

P.S. Some more recent references are added in \cite{Col,Sch,BFKT,Tong,LeSc,KSSY,SoXu,DoEs}.

\ \

Further, the present paper is organized as follows. In the next section we give a synopsis of our approach.
Then we apply this to the group ~$SL(2n,\bbr)$, using the Langlands duality of the subgroup $M$ used in the example.
The cases $n=2,3$ are exposed in separate sections, while the case $n=4$ is briefly summarized  in an Appendix.


\section{Synopsis: Canonical construction of invariant differential operators}

Let $G$ be a semi-simple, non-compact Lie group, and $K$ a
maximal compact subgroup of $G$. Then, we have an {\it Iwasawa
decomposition} $G=KA_0N_0$, where $A_0$ is an Abelian simply
connected vector subgroup of $G$ and $N_0$ is a nilpotent simply
connected subgroup of $G$ preserved by the action of $A_0$.
Furthermore, let $M_0$ be the centralizer of $A_0$ in $K$. Then, the
subgroup $P_0 = M_0 A_0 N_0$ is a {\it minimal parabolic subgroup} of
$G$.  A {\it parabolic subgroup} $P' = M' A' N'$ is any subgroup of $G$
which contains a minimal parabolic subgroup.

Furthermore let $\cg,\ck,\cp,\cm,\ca,\cn$ denote the Lie algebras of $G,K,P,M,A,N$, resp.

Further, for simplicity, we  restrict to {\it maximal parabolic
subgroups} $P=MAN$, i.e.,  $\rank\, A =1$, resp., to {\it maximal parabolic
subalgebras} $\cp = \cm \oplus \ca \oplus \cn$ with $\dim\, \ca=1$.

Let $\nu$ be a (non-unitary) character of $A$, $\nu\in\ca^*$,
parameterized by a real number {\it $d$}, called (for historical reasons) the {\it conformal weight} or
energy.

Furthermore, let $\mu$ fix a discrete series representation
 $D^\mu$ of $M$ on the Hilbert space $V_\mu\,$, or   the
finite-dimensional (non-unitary) representation of $M$ with the same
Casimirs.

 We call the induced
representation $\chi =$ Ind$^G_{P}(\mu\otimes\nu \otimes 1)$ an
 {\it \it elementary representation} of $G$ \cite{DMPPT}. (These are
called {\it generalized principal series representations} (or {\it
limits thereof}) in \cite{Knapp}.)  Their spaces of functions are  \begin{equation}\label{func}
\cc_\chi = \{ \cf \in C^\infty(G,V_\mu) \vr \cf (gman) =
e^{-\nu(H)} \cdot D^\mu(m^{-1})\, \cf (g) \} \end{equation} where $a=
\exp(H)\in A'$,~$H\in\ca'\,$,~$m\in M'$, $n\in N'$. The
representation action is the {\it left regular action}:  \begin{equation}\label{lrega}
(\ct^\chi(g)\cf) (g') = \cf (g^{-1}g'), \quad g,g'\in G.\end{equation}

 An important ingredient in our considerations are the {\it \it
highest/lowest-weight representations} of $\cg^\bac$. These can be
realized as (factor-modules of) Verma modules~$V^\L$ over
 $\cg^\bac$, where $\L\in (\ch^\bac)^*$, $\ch^\bac$ is a Cartan
subalgebra of $\cg^\bac$ and the weight $\L = \L(\chi)$ is determined
uniquely from $\chi$ \cite{VKD1}.

Actually, since our ERs may be induced from finite-dimensional
representations of~$\cm$~(or their limits) the Verma modules are
always reducible. Thus, it is more convenient to use~{\it \it
generalized Verma modules}~$\tV^\L$~such that the role of the
highest/lowest-weight vector $v_0$ is taken by the
(finite-dimensional) space~$V_\mu\,v_0\,$. For the generalized
Verma modules (GVMs) the reducibility is controlled only by the
value of the conformal weight $d$. Relatedly, for the \idos{}, only
the reducibility with regard to non-compact roots is~essential.

 {Another}   main ingredient of our approach is as follows. We group
the (reducible) ERs with the same Casimirs in sets called {\it
multiplets} \cite{VKD1}. The multiplet corresponding to fixed values of the
Casimirs may be depicted as a connected graph, the {\it vertices} of
which correspond to the reducible ERs and the {\it lines (arrows)}
between the vertices correspond to intertwining operators.  The
explicit parameterization of the multiplets and of their ERs is
important in understanding of the situation. The notion of multiplets was introduced in \cite{Dobmul}
and applied to representations of~$SO_o(p,q)$~and~$SU(2,2)$, resp.,
induced from their minimal parabolic subalgebras. Then it was applied to
the conformal superalgebra \cite{DoPemul}, to quantum groups \cite{VKD2}, to
infinite-dimensional (super)algebras \cite{VKD3}.
(For other applications,
we refer to \cite{VKD4}.

In fact, the multiplets contain explicitly all the data necessary to
construct the \idos{}. Actually, the data for each \ido{} consist
of the pair $(\b,m)$, where $\b$ is a (non-compact) positive root
of $\cg^\bac$, $m\in\bbn$, such that the {\it BGG  Verma module
reducibility condition} (for highest-weight modules) is fulfilled:\vspace{-3pt}
\begin{equation}\label{bggr} (\L+\r, \b^\vee ) = m, \quad \b^\vee \equiv 2 \b
/(\b,\b) \ \end{equation}
where $\r$ is half the sum of the positive roots of
 $\cg^\bac$. When the above holds, then the Verma module with shifted
weight $V^{\L-m\b}$ (or $\tV^{\L-m\b}$ for GVM and $\b$
non-compact) is embedded in the Verma module $V^{\L}$ (or
 $\tV^{\L}$). This embedding is realized by a singular vector
 $v_s$ determined by a polynomial $\cp_{m,\b}(\cg^-)$ in the
universal enveloping algebra $(U(\cg_-))\ v_0\,$, and $\cg^-$ is the
subalgebra of $\cg^\bac$ generated by the negative root generators
\cite{Dix}.
 More explicitly \cite{VKD1}, $v^s_{m,\b} = \cp_{m,\b}\, v_0$ (or $v^s_{m,\b} = \cp_{m,\b}\, V_\mu\,v_0$ for GVMs).
   Then,
there exists \cite{VKD1}
 {\bf \ido{}} \begin{equation}\label{invop}  \cd_{m,\b} : \cc_{\chi(\L)}
 \llr \cc_{\chi(\L-m\b)} \end{equation} given explicitly by: \begin{equation}\label{singvv}
 \cd_{m,\b} = \cp_{m,\b}(\widehat{\cg^-})  \end{equation} where
 $\widehat{\cg^-}$ denotes the {\it right action} on the functions
 $\cf$.

\section{Restricted Weyl groups and related notions}

In our exposition below, we shall use the so-called Dynkin labels:  \begin{equation}\label{dynk} m_i
 \equiv (\L+\r,\a^\vee_i) , \quad i=1,\ldots,n, \end{equation} where $\L =
\L(\chi)$, $\r$ is half the sum of the positive roots of
 $\cg^\bac$.

We shall use also   the so-called ~{\it Harish--Chandra parameters}:
\begin{equation}\label{dynhc} m_\b \equiv (\L+\r, \b^\vee )\ ,  \end{equation} where $\b$ is any
positive root of $\cg^\bac$. These parameters are redundant, since
they are expressed in terms of the Dynkin labels; however,   some
statements are best formulated in their terms. (Clearly, both the Dynkin labels and
Harish--Chandra parameters have their origin in the BGG reducibility condition \eqref{bggr}.)

 Next, we recall the action of the Weyl group on highest weights:
\eqn{weeyl}
w_\b(\L) ~\doteq~ \L - (\L+\r,\b^\vee) \b \ee
and thus,
\eqn{wezyl} w_\b(\L) ~=~ \L - m_\b \b \ee
and the shifted weight in \eqref{invop} results by the action of the Weyl group as in \eqref{wezyl}.

Next we mention the important notion of ~{\it restricted Weyl group}. We first need the so-called ~restricted roots.

Let $\D_{\ca'}$ be the {\it restricted root system} of  $(\cg,\ca')$:
\eqnn {rroots} &&\D_{\ca'} ~\doteq~ \{ \l \in \ca'^* ~\vr ~ \l \neq 0 ,\
\cg^\l_{\ca'} \neq 0 \} ~, \nn\\ && \cg^\l_{\ca'} \doteq \{ X \in \cg ~\vr ~
[Y,X] = \l(Y) X ~, ~~\forall Y\in \ca' \} ~. \eea The elements of
$\D_{\ca'}$  are called\ {\it  $\ca'$-restricted roots}.\\
{}[The terminology comes from the fact that things may be arranged so that these roots
are obtained as restriction to $\ca'$ of some roots of the
root system ~$\D$~ of the pair ~$(\cg^\bbc,\ch^\bbc)$.]

  For $\l \in \D_{\ca'}\,$, ~$\cg_{\ca'}^\l$ are called {\it  $\ca'$-restricted root
spaces}, ~dim$_R~\cg_{\ca'}^\l \geq 1$.

Next we introduce some ordering (e.g., the lexicographic one) in $\D_{\ca'}\,$.
Accordingly the latter is split into positive   and
negative   restricted roots: ~$\D_{\ca'} = \D_{\ca'}^{+} \cup \D_{\ca'}^{-}$.

Furthermore, we introduce the simple restricted root system ~$\D^R_{\ca'}\,$,
which is the simple root system of the restricted roots.
Next we introduce the ~{\it restricted Weyl reflections}: for each root ~$\l\in \D_{\ca'}^{+} $~
we define a reflection $s_\l$ in ${\ca'}^*$~:
\eqn{weylres} s_\l (\mu) ~\equiv~ \mu - 2 \frac{(\l,\mu)}{(\l,\l)} \l \ , \quad \mu\in{\ca'}^* \ee
Clearly, ~$s_\l (\l) = -\l\,$, ~$s_\l^2 = $id$_{{\ca'}^*}\,$.

The above reflections generate the ~{\it $\ca'$-restricted Weyl group}~
~$W(\cg,{\ca'})$.

The above may be applied to the case when instead of some ~$\ca'$~ we use an arbitrary subalgebra ~$\ch'$~ of ~$\ch$.


\section{The case of $SL(2n,\bbr)$}

In this paper we treat the case of $G=SL(2n,\bbr)$, $\cg=sl(2n,\bbr)$. We restrict to maximal parabolic subalgebra
\eqnn{parab} \cp ~&=&~ \cm \oplus \ca \oplus \cn \\
\cm &=& sl(n,\bbr) \oplus sl(n,\bbr) \,, ~\dim \ca = 1, ~\dim \cn = n^2 \eea

In the context of relative Langlands duality this case was studied in, e.g.,
 \cite{BSV}, as the subcase of hyperspherical dual pairs
 \cite{JaRa}. There the relation to physics appeared as  arithmetic analog of the electric-magnetic duality of boundary conditions in four-dimensional supersymmetric Yang-Mills theory. This aspect will be recovered for\ $n=2$\ below.

In this section we start with $G=SL(m,\bbr)$, the group of invertible $m \times m $ matrices with real elements and
determinant 1. Then $\cg = sl(m,\bbr)$ and the Cartan involution is given explicitly by:
~$\th X = -\ ^tX$, where $^tX$ is the transpose of $X\in\cg$. Thus, $\ck \cong so(m)$, and is spanned by matrices
(r.l.s. stands for real linear span):
\eqn{maxc}  \ck = {\rm r.l.s.} \{ X_{ij} \equiv e_{ij} - e_{ji} \ , \quad 1\leq i < j \leq m \} \ ,\ee
where $e_{ij}$ are the standard matrices with only nonzero entry (=1) on the $i$-th row and $j$-th column,
$(e_{ij})_{k\ell} = \d_{ik}\d_{j\ell}\,$. (Note that $\cg$ does not have discrete series representations
if $m>2$.)

Further, the complementary space $\cp$ is given by:
\eqnn{slp} \cp &=& {\rm r.l.s.} \{ Y_{ij} \equiv e_{ij} + e_{ji} \ , ~~ 1\leq i < j \leq m \ ,\\
&&  H_j \equiv e_{jj} - e_{j+1,j+1} \ , \quad 1\leq j \leq m-1 \}\ .\eea
The split rank is ~$r = m-1$, and from \eqref{slp} it is obvious that in this setting one has:
\eqn{aaa}\ch = {\rm r.l.s.} \{ H_j   \ , \quad 1\leq j \leq n-1=r \}\ .\ee

The simple root vectors are given explicitly by:
\eqn{simple} X^+_j \doteq e_{j,j+1} \ , ~~~ X^-_j \doteq e_{j+1,j} \ , ~~~ 1\leq j \leq m-1 \ .\ee
Note that matters are arranged so that
\eqn{slsub} [X^+_j,X^-_j] = H_j\ , \quad [H_j,X^\pm_j] = \pm 2X^\pm_j\ ,\ee
and further we shall denote by ~$sl(2,\bbr)_j$~ the ~$sl(2,\bbr)$~
subalgebra of $\cg$ spanned  by $X^\pm_j\,, H_j\,$.

In our case of consideration ~$m=2n$ we have
\eqn{facm}  \cm ~=~ sl(n) \oplus sl(n) \ee
and we use representations of $\cm$ indexed as follows:
\eqn{repm} \hat{\cm} = (m_1,\ldots,m_{n-1} \ ;\ m_{n+1}, \ldots, m_{2n-1}) \ee
When all $m_j$ are natural numbers   $\hat{\cm}$ indexes the unitary finite-dimensional irreps of $\cm$.


\subsection{$sl(4)$}

In the case of $sl(4)$ the 
  parabolic $\cm$ factor is:
\eqn{su22m} \cm_4 ~=~ sl(2) \oplus sl(2) \ee
the representations being indexed by the numbers $m_1,m_3$

Relatedly the representations of $\cg$ are indexed by:
\eqn{su22g} \chi_4 ~=~ [m_1,m_2,m_3] \ee
It is well-known that when all $m_j$ are natural numbers then $\chi_4$ exhausts the finite-dimensional representations of $\cg$. 
Each representation $\chi_4$ is part of 24-member multiplet naturally corresponding to the 24 elements of the Weyl group of $sl(4)$.
When we consider   induction from $\cm_4$ then we have six-member multiplets (sextets) parametrized as follows:
\eqnn{tabsu22}
\chi^- ~&=&~ \{\, m_1\,, m_2\,, m_3\,
\}  , \\
\chi'^- ~&=&~ \{\, m_{12}\,, -m_{2}\,, m_{23}
\}  , \ \L'^- = \L^- - m_2\a_2 \nn\\
 \chi''^- ~&=&~ \{\, m_{2}\,, -m_{12}\,, m_{13}
\}  , \ \L''^-  = \L'^- - m_1\a_{12} \nn\\
\chi''^+ ~&=&~ \{\, m_{13}\,, -m_{23}\,, m_{2}
\}  , \ \L''^+ = \L'^- - m_3\a_{23} \nn\\
 \chi'^+ ~&=&~ \{\, m_{23}\,, -m_{13}\,, m_{12}
\}  , \ \L'^+ = \L''^- - m_3\a_{23} = \L''^+ - m_1\a_{12} \nn\\
 \chi^+ ~&=&~ \{\, m_{3}\,, -m_{13}\,, m_{1}
\}  , \ \L^+ = \L'^+ - m_2\a_{2} \nn\eea
where  $m_{12}\equiv m_1+m_2$, $m_{23}\equiv m_2+m_3$, $m_{13}\equiv m_1+m_2+m_3$.
Note that the $\pm$ pairs are related by Knapp=Stein \cite{KnSt} integral intertwining operators ~$G^\pm$~ so that
the operators $G^+$ act from $\chi^-$ to $\chi^+$, while $G^-$ act from $\chi^+$ to $\chi^-$, etc.

Thus, the ~{\it Knapp-Stein duality}~ is a manifestation of the ~{\it Langlands duality}.

We recall  that the number ~$N_M$~ of ERs in a multiplet corresponding to induction from a parabolic given by \cite{VKD1}:
\eqn{redun} N_M ~=~  {\vert W(\cg,\ch) \vert \over \vert W(\cm,\ch_m) \vert} \ee
which in our case ($\cm=\cm_4$) gives:
\eqn{redu22} N_M ~=~  {24 \over 4} ~=~ 6. \ee
what we have obtained.

An alternative parametrization stressing the duality is given as follows:
\eqnn{tabsu222}
\chi^\pm  ~&=&~ \{\, (m_1\,;   m_3)^\pm\,;\ c ~=~ \pm (m_2 +\ha (m_1+m_3)) \
\}  , \nn\\
\chi'^\pm ~&=&~ \{\, (m_{12}\,,   m_{23})^\pm \,;\ c ~=~ \pm \ha (m_1+m_3)
\}  , \nn\\
 \chi''^\pm ~&=&~ \{\,( m_{2}\,,   m_{13})^\pm \,;\ c ~=~ \pm \ha (m_1-m_3)
\}  , \
  \nn\eea
  where ~$(p;q)^+ = (q;p)$, ~$(p;q)^- = (p;q)$,

The irreducible subrepresentations ~$\ce$~ of ~$\chi^-$~ are finite-dimensional, exhausting all finite-dimensional (non-unitary) representations of ~$sl(4)$, and of all real forms. 

Note also that the dimensions of the $\pm$ inducing pair of $\cm$ are the same, namely, $m_1m_3$ for $\chi^\pm$,
$m_{12}m_{23}$ for $\chi'^\pm$, $m_{2}m_{13}$ for $\chi''^\pm$.

Finally, we use the simplest case ~$m_1=m_2=m_3=1$~ to exhibit the electro-magnetic duality which has transparent
physical meaning for the conformal real form $su(2,2)$.
The multiplet is depicted on Fig. 1. (Complete treatment may be found in \cite{VKD1}.)





\fig{}{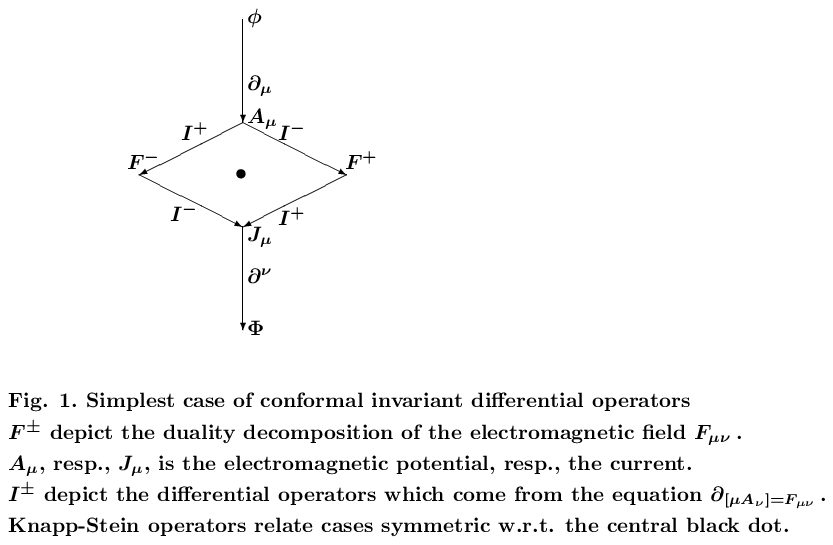}{14cm}

Multiplets containing the finite-dimensional subrepresentations are called {\it main multiplets}.
 The other multiplets are called reduced multiplets.
These contain inducing finite-dimensional representations of $\cm$.

In the case at hand there are three such cases so that each is a doublet (containing two ERs).
 Explicitly the three cases are:
\eqna{tabsu221}
 \chi_1^\pm ~&=&~ \{\, (m_{2}\,,   m_{23})^\pm \,;\ c ~=~ \pm \ha m_3\} ,\,
  \L_1^+ = \L_1^- - m_3\a_{23}\\
  \chi_2^\pm  ~&=&~ \{\, (m_1\,;   m_3)^\pm\,;\ c^\pm ~=~ \pm \ha (m_1+m_3)
\}  , \\
 \chi_3^\pm ~&=&~ \{\, (m_{12}\,,   m_{2})^\pm \,;\ c ~=~ \pm \ha m_1\} ,\,
 \L_3^+ = \L_3^- - m_1\a_{12}\eena
   Note that here the invariant operators
are deformations of the   Knapp-Stein integral operators from the sextet picture.
Thus, those from ~$\chi^+$~ to ~$\chi^-$~ are still integral operators, while those
from ~$\chi^-$~ to ~$\chi^+$~ are differential operators via
degeneration of the Knapp-Stein integral operators. Yet in the first and third case
 these are differential operators inherited from the sextets,
 only the  operators in (\ref{tabsu221}b) from  ~$\chi_2^-$~ to ~$\chi_2^+$~  are
obtained due to genuine degeneration of the Knapp-Stein integral operators \cite{VKD1}.
This is the standard degeneration of the two-point function-kernel which at the reducibility points is a generalized
function with regularization turning it into delta-function (cf. Gelfand et al \cite{GeV}).
Finally, we  add that in the case ~$m_1=m_3=n$~  the operators (\ref{tabsu221}b)
become a degree of the d'Alembert operator:
\eqn{degdala}  \cd_{n,n} ~=~ const\, \Box^{c^+} ~=~ const\, \Box^n \ee


\section{$sl(6)$}

 Here we take up the case $sl(6)$   with parabolic $\cm$ factor
 \eqn{cmm33} \cm_5 ~=~ sl(3,\bbr)\oplus sl(3,\bbr) ~=~ \cm_{2L} \oplus  \cm_{2R} \ee

We start with representations of $sl(6,\bbr)$ indexed by five numbers:
\eqn{fivem} \chi ~=~ \{ m_1\,, m_2\,,m_3\,,m_4\,,m_5\}\,, \ee
so that ~$m_1,m_2$~ index the representations of  $\cm_{2L}$,
~$m_4,m_5$~ index the representations of  $\cm_{2R}$, while ~$m_3$~
indexes the representations of the dilatation subalgebra $\ca$.

When all $m_j$ are positive integers we use   formula \eqref{redun} so
we have a multiplet of 20   members  since:
\eqn{redu66}
 N_M ~=~  {\vert W(\cg,\ch) \vert \over \vert W(\cm_5,\ch_5) \vert} ~=~  {6! \over (3!)^2} ~=~ 20. \ee

Explicitly, the signatures are:
\eqnn{tabl6}
&&\chi_1  ~=~ \{\, m_1\,, m_2\,, m_3\,,m_4\,, m_5\,
\}  , \\
&&\chi_2 ~=~ \{\, m_1\,, m_{23}\,, -m_3\,,m_{34}\,, m_5\,
\}  , \ \L_2 = \L_1 - m_3\a_3 \nn\\
 &&\chi_3 ~=~ \{\, m_{12}\,, m_{3}\,, -m_{23}\,,m_{24}\,, m_5\,
\}  , \ \L_3 = \L_2 - m_2\a_{23} \nn\\
&&\chi_4 ~=~ \{\, m_1\,, m_{24}\,, -m_{34}\,,m_{3}\,, m_{45}\,
\}  , \ \L_4 = \L_2 - m_4\a_{34} \nn\\
 &&\chi_5 ~=~ \{\, m_{2}\,, m_{3}\,, -m_{13}\,,m_{14}\,, m_5\,
\}  , \ \L_5 = \L_3 - m_1\a_{13} \nn\\
 &&\chi_6 ~=~ \{\, m_{12}\,, m_{34}\,, -m_{24}\,,m_{23}\,, m_{45}\,
\}  , \ \L_6 = \L_3 - m_4\a_{34} \nn\\
&&\chi_7 ~=~ \{\, m_1\,, m_{25}\,, -m_{35}\,,m_{3}\,, m_{4}\,
\}  , \ \L_7 = \L_4 - m_5\a_{35} \nn\\
&&\chi_8 ~=~ \{\, m_{2}\,, m_{34}\,, -m_{14}\,,m_{13}\,, m_{45}\,
\}  , \ \L_8 = \L_6 - m_1\a_{13} \nn\\
&&\chi_9 ~=~ \{\, m_{13}\,, m_{4}\,, -m_{24}\,,m_{2}\,, m_{35}\,
\}  , \ \L_9 = \L_6 - m_3\a_{24} \nn\\
&&\chi_{10} ~=~ \{\, m_{12}\,, m_{35}\,, -m_{25}\,,m_{23}\,, m_{4}\,
\}  , \ \L_{10} = \L_6 - m_5\a_{35} \nn\\
&&\chi_{11} ~=~ \{\, m_{23}\,, m_{4}\,, -m_{14}\,,m_{12}\,, m_{35}\,
\}  , \ \L_{11} = \L_8 - m_3\a_{24} \nn\\
&&\chi_{12} ~=~ \{\, m_{2}\,, m_{35}\,, -m_{15}\,,m_{13}\,, m_{4}\,
\}  , \ \L_{12} = \L_8 - m_5\a_{35} \nn\\
&&\chi_{13} ~=~ \{\, m_{13}\,, m_{45}\,, -m_{25}\,,m_{2}\,, m_{34}\,
\}  , \ \L_{13} = \L_{10} - m_3\a_{24} \nn\\
&&\chi_{14} ~=~ \{\, m_{3}\,, m_{4}\,, -m_{14}\,,m_{1}\,, m_{25}\,
\}  , \ \L_{3} = \L_{11} - m_2\a_{14} \nn\\
&&\chi_{15} ~=~ \{\, m_{23}\,, m_{45}\,, -m_{15}\,,m_{12}\,, m_{34}\,
\}  , \ \L_{15} = \L_{11} - m_5\a_{35} \nn\\
&&\chi_{16} ~=~ \{\, m_{14}\,, m_{5}\,, -m_{25}\,,m_{2}\,, m_{3}\,
\}  , \ \L_{13} = \L_{10} - m_4\a_{25} \nn\\
&&\chi_{17} ~=~ \{\, m_{3}\,, m_{45}\,, -m_{15}\,,m_{1}\,, m_{24}\,
\}  , \ \L_{17} = \L_{15} - m_2\a_{14} \nn\\
&&\chi_{18} ~=~ \{\, m_{24}\,, m_{5}\,, -m_{15}\,,m_{12}\,, m_{3}\,
\}  , \ \L_{18} = \L_{15} - m_4\a_{25} \nn\\
&&\chi_{19} ~=~ \{\, m_{34}\,, m_{5}\,, -m_{15}\,,m_{1}\,, m_{23}\,
\}  , \ \L_{19} = \L_{17} - m_4\a_{25} \nn\\
&&\chi_{20} ~=~ \{\, m_{4}\,, m_{5}\,, -m_{15}\,,m_{1}\,, m_{2}\,
\}  , \ \L_{20} = \L_{19} - m_3\a_{15} \nn
\eea
Note that we have indicated some embeddings between Verma modules but not all in order not to clutter the formulae. The full picture is seen on Fig. 2.


\fig{}{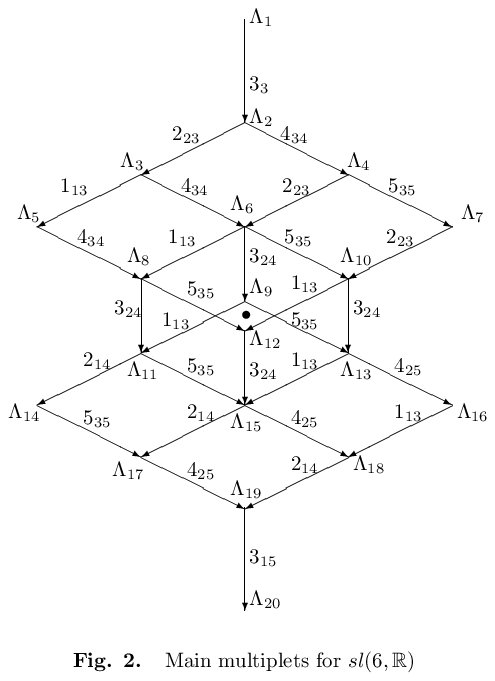}{20cm}

We quickly observe that the representations ~$\chi_{n}$~ and ~$\chi_{21-n}$~ are Langlands duals related by  Knapp-Stein operators. More explicitly, this duality is given by the following presentation of the same
multiplet:
\eqnn{tabl6pm}
&&\chi^\pm_1  ~=~ \{\, (m_1\,, m_2\,;\, m_4\,, m_5)^\pm\,;\, c=\pm (m_3 + \ha m_{12,45})\,
\}  , \\
&&\chi^\pm_2 ~=~ \{\, (m_1\,, m_{23}\,;\, m_{34}\,, m_5)^\pm\,;\,  c= \pm\ha m_{12,45}\,
\}   ,\nn\\ 
 &&\chi^\pm_3 ~=~ \{\, (m_{12}\,, m_{3}\,;\,m_{24}\,, m_5)^\pm\,;\,  c= \pm\ha m_{1,45}\,
\}   ,\nn\\ 
&&\chi^\pm_4 ~=~ \{\, (m_1\,, m_{24}\,; \,m_{3}\,, m_{45})^\pm\,;\,  c= \pm\ha m_{12,5}\,
\}   ,\nn\\ 
 &&\chi^\pm_5 ~=~ \{\, (m_{2}\,, m_{3}\,; \,m_{14}\,, m_5)^\pm\,;\,  c= \pm\ha (m_{45}-m_1)\,
\}   ,\nn\\ 
 &&\chi^\pm_6 ~=~ \{\, (m_{12}\,, m_{34}\,; \,m_{23}\,, m_{45})^\pm\,;\,  c= \pm\ha (m_1+m_{5}) \,
\}   ,\nn\\ 
&&\chi^\pm_7 ~=~ \{\, (m_1\,, m_{25}\, ;\,m_{3}\,, m_{4} )^\pm\,;\,  c= \pm\ha (m_{12}-m_{5})\,
\}   ,\nn\\ 
&&\chi^\pm_8 ~=~ \{\,( m_{2}\,, m_{34}\, ; \,m_{13}\,, m_{45})^\pm\,;\,  c= \pm\ha (m_{5}-m_1)\, \,
\}   ,\nn\\ 
&&\chi^\pm_9 ~=~ \{\, (m_{13}\,, m_{4}\,;\,m_{2}\,, m_{35})^\pm\,;\,  c= \pm\ha (m_1+m_{5})\,
\}   ,\nn\\ 
&&\chi^\pm_{10} ~=~ \{\, (m_{12}\,, m_{35}\,;\, m_{23}\,, m_{4})^\pm\,;\,  c= \pm\ha (m_{1}-m_{5})\,
\}   ,\nn 
\eea
where ~$(p,q;r,s)^+ \equiv (r,s;p,q)$, ~
~$(p,q;r,s)^- \equiv (p,q;r,s)$, and
 the inducing number of the dilatation subalgebra ~$\ca$~ is replaced by the conformal factor ~$c$.
 Clearly,  ~$\chi^-_n = \chi_n$, ~$\chi^+_n = \chi_{21-n}$~ ~ for ~$1\leq n \leq 10$.

\subsection*{Reduced multiplets}

 Here  we just list the reduced multiplets:
\eqna{tabl6pm1}
  &1&\chi'^\pm_3 ~=~ \{\, (m_{2}\,, m_{3}\,;\,m_{24}\,, m_5)^\pm\,;\,  c= \pm\ha m_{45}\,
\}   ,\\ 
   &&\chi'^\pm_6 ~=~ \{\, (m_{2}\,, m_{34}\,; \,m_{23}\,, m_{45})^\pm\,;\,  c= \pm\ha m_{5} \,
\}   ,\nn\\ 
  &&\chi'^\pm_9 ~=~ \{\, (m_{23}\,, m_{4}\,;\,m_{2}\,, m_{35})^\pm\,;\,  c= \pm\ha m_{5}\,
\}   ,\nn\\ 
  &13&\chi''^\pm_6 ~=~ \{\, (m_{2}\,, m_{4}\,; \,m_{2}\,, m_{45})^\pm\,;\,  c= \pm\ha m_{5} \,
\}   ,\nn\\ 
 &14&\chi''^\pm_3 ~=~ \{\, (m_{2}\,, m_{3}\,;\,m_{23}\,, m_5)^\pm\,;\,  c= \pm\ha m_{5}\,
\}   ,\nn\\ 
&15&\chi''^\pm_9 ~=~ \{\, (m_{23}\,, m_{4}\,; \,m_{2}\,, m_{34})^\pm\,;\,  c=  0 \,
\}   ,\\ 
 &135&\chi_6 ~=~ \{\, (m_{2}\,, m_{4}\,; \,m_{2}\,, m_{4})\,;\,  c=  0 \,
\}   ,   
  \eena
       \eqnn{tabl6pm2}
 &2&\chi'^\pm_2 ~=~ \{\, (m_1\,, m_{3}\,;\, m_{34}\,, m_5)^\pm\,;\,  c= \pm\ha m_{1,45}\,
\}   ,\\ 
 &&\chi'^\pm_4 ~=~ \{\, (m_1\,, m_{34}\,; \,m_{3}\,, m_{45})^\pm\,;\,  c= \pm\ha m_{1,5}\,
\}   ,\nn\\ 
 &&\chi'^\pm_7 ~=~ \{\, (m_1\,, m_{35}\, ;\,m_{3}\,, m_{4} )^\pm\,;\,  c= \pm\ha (m_{1}-m_{5})\,
\}   ,\nn\\ 
  &24&\chi''^\pm_2 ~=~ \{\, (m_{1}\,, m_{3}\,;\,m_{3}\,, m_5)^\pm\,;\,  c= \pm\ha m_{1,5}\,
\}   ,\nn\\ 
 &25&\chi''^\pm_4 ~=~ \{\, (m_{1}\,, m_{34}\,; \,m_{3}\,, m_{4})^\pm\,;\,  c= \pm\ha m_1 \,
\}   ,\nn  
\eea
\eqnn{tabl6pm3}
  &3&\chi'^\pm_1 ~=~ \{\, (m_1\,, m_{2}\,;\, m_{4}\,, m_5)^\pm\,;\,  c= \pm\ha m_{12,45}\,
\}   ,\\ 
  &&\chi'^\pm_6 ~=~ \{\, (m_{12}\,, m_{4}\,; \,m_{2}\,, m_{45})^\pm\,;\,  c= \pm\ha m_{1,5} \,
\}   ,\nn\\ 
 &&\chi'^\pm_8 ~=~ \{\,( m_{2}\,, m_{4}\, ; \,m_{12}\,, m_{45})^\pm\,;\,  c= \pm\ha (m_{5}-m_1)\, \,
\}   ,\nn 
   \eea
  Note that the numbers on the left indicate which representation numbers are missing  in the displayed signatures.\\
Further, note that the ~$\pm$~ pairs are Knapp-Stein pairs, except the case
(\ref{tabl6pm1}b) where the operator is just  a flip of the finite-dimensional inducing irreps. Note also that the case (\ref{tabl6pm1}c)  is a singlet.\\  Note that we do not display  reduced multiplets with missing labels $m_4$ and $m_5$ since due to duality they are equivalent to multiplets with missing labels $m_2$, $m_1$, resp.

\section*{Outlook} Obviously, there many more cases of Langlands duals to which our approach  can be applied. This will be done in some future papers.
Note also that the cases ~$sl(2n)$~ may be treated as ~$su(n,n)$~ as displayed above for $n=2$.


\section*{Appendix: sl(8)}

Here we consider briefly  the case $sl(8)$   with parabolic $\cm$ factor
 \eqn{cmm44} \cm_7 ~=~ sl(4,\bbr)\oplus sl(4,\bbr)  ~=~ \cm_{2L} \oplus  \cm_{2R} \ee

Analogously to the previously considered cases the representations of $sl(9,\bbr)$ indexed by seven numbers:
\eqn{sevem} \chi ~=~ \{ m_1\,, m_2\,,m_3\,,m_4\,,m_5\,,m_6,\,,m_7 \}\,, \ee
so that ~$m_1,m_2,m_3$~ index the representations of  $\cm_{2L}$,
~$m_5,m_6,m_7$~ index the representations of  $\cm_{2R}$, and ~$m_4$~
indexes the representations of the dilatation subalgebra $\ca$.

When all $m_j$ are positive integers we again use the formula \eqref{redun} so
we have a multiplet of 70   members  since:
\eqn{redu4}
 N_M ~=~  {\vert W(\cg,\ch) \vert \over \vert W(\cm_7,\ch_7) \vert} ~=~  {8! \over (4!)^2} ~=~ 70. \ee

For shortness we give only the Knapp-Stein dual signatures:
\eqnn{tabsl8}
&&\chi^\pm_1  ~=~ \{\, (m_1\,, m_2\,, m_3\,,m_4\,, m_5\,, m_6\,,m_7)^\pm  , \ c^\pm ~=~ \pm (m_4 + \ha m_{13,57})\,
\}  , \nn\\
&&\chi^\pm_2  ~=~ \{\,(m_1\,, m_2\,, m_{34}\,,-m_4\,, m_{45}\,, m_6\,,m_7)^\pm \,
   ,\  \ c^\pm ~=~ \pm ( \ha m_{17})\} \\ 
&&\chi^\pm_3  ~=~ \{\,(m_1\,, m_{23}\,, m_4\,,-m_{34}\,, m_{35}\,, m_6\,,m_7)^\pm\ , \ c^\pm ~=~ \pm ( \ha m_{12,57})\, \}\nn\\ 
&&\chi^\pm_4  ~=~ \{\,(m_1\,, m_2\,, m_{35}\,,-m_{45}\,, m_4\,, m_{56}\,,m_7)^\pm\ , \ c^\pm ~=~ \pm ( \ha m_{13,67})\, \}\nn\\ 
&&\chi^\pm_5  ~=~ \{\,(m_{12}\,, m_{3}\,, m_4\,,-m_{24}\,, m_{25}\,, m_6\,,m_7)^\pm\ , \ c^\pm ~=~ \pm ( \ha m_{1,57})\,  \}\nn\\ 
&&\chi^\pm_6  ~=~ \{\,(m_{1}\,, m_{23}\,, m_{45}\,,-m_{35}\,, m_{34}\,, m_{56}\,,m_7)^\pm\ , \ c^\pm ~=~ \pm ( \ha m_{12,67})\, \}\nn\\ 
&&\chi^\pm_7  ~=~ \{\,(m_{1}\,, m_{2}\,, m_{36}\,,-m_{46}\,, m_{4}\,, m_{5}\,, m_{67})^\pm\ , \ c^\pm ~=~ \pm ( \ha m_{13,7})\, \}\nn\\ 
&&\chi^\pm_8  ~=~ \{\,(m_{2}\,, m_{3}\,, m_4\,,-m_{14}\,(m_{15}\,, m_{56}\,,m_7)^\pm\ , \ c^\pm ~=~ \pm ( \ha m_{-1,57})\, \}\nn\\ 
&&\chi^\pm_9  ~=~ \{\,(m_{12}\,, m_{3}\,, m_{45}\,,-m_{25}\,(m_{24}\,, m_{56}\,,m_7)^\pm\ , \ c^\pm ~=~ \pm ( \ha m_{1,67})\, \}\nn\\ 
&&\chi^\pm_{10}  ~=~ \{\,(m_{1}\,, m_{24}\,, m_{5}\,,-m_{35}\,(m_{3}\,, m_{46}\,,m_7)^\pm\ , \ c^\pm ~=~ \pm ( \ha m_{12,67})\,  \}\nn\\ 
&&\chi^\pm_{11}  ~=~ \{\,(m_{1}\,, m_{23}\,, m_{46}\,,-m_{36}\,(m_{34}\,, m_{5}\,,m_{67})^\pm\ , \ c^\pm ~=~ \pm ( \ha m_{12,7})\, \}\nn\\ 
&&\chi^\pm_{12}  ~=~ \{\,(m_{1}\,, m_{2}\,, m_{37}\,,-m_{47}\,, m_{4}\,, m_{5}\,, m_{6})^\pm\ , \ c^\pm ~=~ \pm ( \ha m_{13,-7})\, \}\nn\\ 
&&\chi^\pm_{13}  ~=~ \{\,(m_{2}\,, m_{3}\,, m_{45}\,,-m_{15}\,(m_{14}\,, m_{56}\,,m_{7})^\pm\ , \ c^\pm ~=~ \pm ( \ha m_{-1,67})\, \}\nn\\ 
&&\chi^\pm_{14}  ~=~ \{\,(m_{12}\,, m_{34}\,, m_{5}\,,-m_{25}\,(m_{23}\,, m_{46}\,,m_7)^\pm\ , \ c^\pm ~=~ \pm ( \ha m_{1,67})\,  \}\nn\\ 
&&\chi^\pm_{15}  ~=~ \{\,(m_{12}\,, m_{3}\,, m_{46}\,,-m_{26}\,(m_{24}\,, m_{5}\,,m_{67})^\pm\ , \ c^\pm ~=~  \pm ( \ha m_{1,7})\, \}\nn\\ 
&&\chi^\pm_{16}  ~=~ \{\,(m_{1}\,, m_{24}\,, m_{56}\,,-m_{36}\,(m_{3}\,, m_{45}\,,m_{67})^\pm\ , \ c^\pm ~=~ \pm ( \ha m_{12,7})\, \}\nn\\ 
&&\chi^\pm_{17}  ~=~ \{\,(m_{1}\,, m_{23}\,, m_{47}\,,-m_{37}\,, m_{34}\,, m_{5}\,, m_{6})^\pm\ , \ c^\pm ~=~ \pm ( \ha m_{12,-7})\, \}\nn\\ 
&&\chi^\pm_{18}  ~=~ \{\,(m_{2}\,, m_{3}\,, m_{46}\,,-m_{16}\,, m_{14}\,, m_{5}\,,m_{67})^\pm\ , \ c^\pm ~=~ \pm ( \ha m_{-1,7})\, \}\nn\\
&&\chi^\pm_{19}  ~=~ \{\,(m_{2}\,, m_{34}\,, m_{5}\,,-m_{15}\,, m_{13}\,, m_{46}\,,m_{7})^\pm\ , \ c^\pm ~=~ \pm ( \ha m_{-1,67})\, \}\nn\\ %
&&\chi^\pm_{20}  ~=~ \{\,(m_{13}\,, m_{4}\,, m_{5}\,,-m_{25}\,, m_{2}\,, m_{36}\,,m_{7})^\pm\ , \ c^\pm ~=~ \pm ( \ha m_{1,67})\, \}\nn\\
&&\chi^\pm_{21}  ~=~ \{\,(m_{12}\,, m_{34}\,, m_{56}\,,-m_{26}\,, m_{23}\,, m_{45}\,,m_{67})^\pm\ ,
\ c^\pm ~=~ \pm ( \ha m_{1,7})\,  \}\nn\\
&&\chi^\pm_{22}  ~=~ \{\,(m_{1}\,, m_{25}\,, m_{6}\,,-m_{36}\,, m_{3}\,, m_{4}\,,m_{57})^\pm\ , \
c^\pm ~=~ \pm ( \ha m_{12,7})\,  \}\nn\\
&&\chi^\pm_{23}  ~=~ \{\,(m_{1}\,, m_{24}\,, m_{57}\,,-m_{37}\,, m_{3}\,, m_{45}\,,m_{6})^\pm\ , \
c^\pm ~=~ \pm ( \ha m_{12,-7})\, \}\nn\\
&&\chi^\pm_{24}  ~=~ \{\,(m_{12}\,, m_{3}\,, m_{47}\,,-m_{27}\,, m_{24}\,, m_{5}\,, m_{6})^\pm\ , \
c^\pm ~=~ \pm ( \ha m_{1,-7})\, \}\nn\\
&&\chi^\pm_{25}  ~=~ \{\,(m_{2}\,, m_{3}\,, m_{47}\,,-m_{17}\,, m_{14}\,, m_{5}\,,m_{6})^\pm\ , \
c^\pm ~=~ \mp ( \ha m_{1,7})\,  \}\nn\\
&&\chi^\pm_{26}  ~=~ \{\,(m_{2}\,, m_{34}\,, m_{56}\,,-m_{16}\,, m_{13}\,, m_{45}\,,m_{67})^\pm\ , \
c^\pm ~=~ \pm ( \ha m_{-1,7})\, \}\nn\\
&&\chi^\pm_{27}  ~=~ \{\,(m_{23}\,, m_{4}\,, m_{5}\,,-m_{15}\,, m_{12}\,, m_{36}\,,m_{7})^\pm\ , \
c^\pm ~=~ \pm ( \ha m_{-1,67})\, \}\nn\\ 
&&\chi^\pm_{28}  ~=~ \{\,(m_{13}\,, m_{4}\,, m_{56}\,,-m_{26}\,, m_{2}\,, m_{35}\,,m_{67})^\pm\ , \
c^\pm ~=~ \pm ( \ha m_{1,7})\, \}\nn\\ 
&&\chi^\pm_{29}  ~=~ \{\,(m_{12}\,, m_{35}\,, m_{6}\,,-m_{26}\,, m_{23}\,, m_{4}\,,m_{57})^\pm\ , \
c^\pm ~=~ \pm ( \ha m_{1,7})\, \}\nn\\ 
&&\chi^\pm_{30}  ~=~ \{\,(m_{1}\,, m_{25}\,, m_{67}\,,-m_{37}\,, m_{3}\,, m_{4}\,,m_{56})^\pm\ , \
c^\pm ~=~ \pm ( \ha m_{12,-7})\, \}\nn\\ 
&&\chi^\pm_{31}  ~=~ \{\,(m_{12}\,, m_{34}\,, m_{57}\,,-m_{27}\,, m_{23}\,, m_{45}\,,m_{6})^\pm\ , \
c^\pm ~=~ \pm ( \ha m_{1,-7})\, \}\nn\\ 
 &&\chi^\pm_{32}  ~=~ \{\,(m_{2}\,, m_{34}\,, m_{57}\,,-m_{17}\,, m_{13}\,, m_{45}\,,m_{6})^\pm\ , \
 c^\pm ~=~ \mp ( \ha m_{1,7})\, \}\nn\\ 
&&\chi^\pm_{33}  ~=~ \{\,(m_{2}\,, m_{35}\,, m_{6}\,,-m_{16}\,, m_{13}\,, m_{4}\,,m_{57})^\pm\ , \
c^\pm ~=~ \pm ( \ha m_{-1,7})\, \}\nn\\ 
&&\chi^\pm_{34}  ~=~ \{\,(m_{23}\,, m_{4}\,, m_{56}\,,-m_{16}\,, m_{12}\,, m_{35}\,,m_{67})^\pm\ , \
c^\pm ~=~ \pm ( \ha m_{-1,7})\, \}\nn\\ 
&&\chi^\pm_{35}  ~=~ \{\,(m_{3}\,, m_{4}\,, m_{5}\,,-m_{15}\,, m_{1}\,, m_{26}\,,m_{7})^\pm
\ , \ c^\pm ~=~ \pm ( \ha m_{-12,67})\, \}\nn\\ 
\nn\eea
where ~$(p,q,u;r,s,v)^+ \equiv (r,s,v;p,q,u)$, ~
~$(p,q,u;r,s,v)^- \equiv (p,q,u;r,s,v)$.

\end{document}